\documentclass[10pt]{article}
\usepackage{amsmath}
\usepackage{amssymb}
\usepackage{amsthm}
\usepackage{amsfonts}
\usepackage{verbatim}
\usepackage{url}
\theoremstyle{plain}
\newtheorem{theorem}{Theorem}[section]

\newtheorem{lemma}[theorem]{Lemma}
\newtheorem{corollary}[theorem]{Corollary}
\newtheorem{proposition}[theorem]{Proposition}

\theoremstyle{definition}
\newtheorem{remark}{Remark}[section]
\newtheorem{example}{Example}[section]
\newtheorem*{notation}{Notation}

\theoremstyle{remark}
\newtheorem*{note}{Note}

\newcommand{\cA}{{\mathcal A}}

\newcommand{\cD}{{\mathcal D}}

\newcommand{\cF}{{\mathcal F}}
\newcommand{\cR}{{\mathcal R}}
\newcommand{\NN}{{\mathbb N}}
\newcommand{\cAstar}{\mathcal{A}^*}
\newcommand{\cAplus}{\mathcal{A}^+}
\newcommand{\empt}{\varepsilon}
\newcommand{\cAw}{\mathcal{A}^{\omega}}
\newcommand{\cAinf}{\mathcal{A}^{\infty}}

\newcommand{\rev}{\widetilde}

\newcommand{\ee}{\'{e}\'{e}}
\newcommand{\pref}{\subseteq_p}
\newcommand{\suff}{\subseteq_s}

\newcommand{\ov}{\overline}
\newcommand{\bs}{\mathbf{s}}
\newcommand{\bx}{\mathbf{x}}
\newcommand{\bt}{\mathbf{t}}
\newcommand{\bff}{\mathbf{f}}
\newcommand{\bz}{\mathbf{z}}

\numberwithin{equation}{section}

\author{Amy Glen\footnotemark[1]}
\title{Occurrences of palindromes in characteristic Sturmian words}
\date{June 24, 2005}
\begin{document}
\normalsize
\begin{comment}
%-------------------------Page 1---------------------------------------------
\begin{titlepage}
\begin{flushleft} \vspace{2cm}
\hrule \vspace{3.5cm}

{\huge{\bf Occurrences of palindromes in}} \\ \vspace{0.5cm}
{\huge{\bf characteristic Sturmian words}}

\vspace*{3cm}

{\LARGE Amy Glen}

\vspace{1.5cm}

{\large School of Mathematical Sciences,} 
\\ \vspace{0.3cm}
{\large Discipline of Pure Mathematics,}
\\ \vspace{0.3cm}
{\large University of Adelaide, South Australia, 5005}

\vspace{1.5cm}

{\large E-mail: \texttt{amy.glen@adelaide.edu.au}} 
\\ \vspace{0.3cm}
{\large Phone: $+61~8~8303~3026$}
\\ \vspace{0.3cm}
{\large Fax: $+61~8~8303~3696$}

\vspace{1.5cm}

{\Large June 2005}

\vspace{6cm} \hrule

\end{flushleft}
\end{titlepage}
\newpage
\end{comment}
%--------------------Beginning of article-----------------------------------
\normalsize \maketitle \footnotetext[1]{E-mail:
\texttt{amy.glen@adelaide.edu.au}} 
\begin{center}
\vspace{-0.7cm} School of Mathematical Sciences, 
Discipline of Pure Mathematics, University of Adelaide, \\ 
South Australia, Australia, 5005
\end{center}
\vspace{0.5cm} \hrule
\begin{abstract}
This paper is concerned with palindromes occurring in
characteristic Sturmian words $c_\alpha$ of slope $\alpha$, 
where $\alpha \in (0,1)$ is an irrational. 
As $c_\alpha$ is a uniformly recurrent infinite word, any (palindromic) factor of $c_\alpha$  
occurs infinitely many times in $c_\alpha$ with bounded gaps. 
Our aim is to completely describe where palindromes occur in $c_\alpha$. 
In particular, given any palindromic factor $u$ of
$c_\alpha$, we shall establish a decomposition of $c_\alpha$ with
respect to the occurrences of $u$. Such a
decomposition shows precisely where $u$ occurs in $c_\alpha$, and
this is directly related to the continued fraction expansion of
$\alpha$.
\vspace{0.1cm} 
\vspace{0.2cm}
\end{abstract}
\hrule

\section{Introduction}
The fascinating family of Sturmian words consists of all aperiodic
infinite words having exactly $n+1$ distinct factors of length $n$
for each $n \in \NN$. Such words have many applications in
various fields of mathematics, such as symbolic dynamics, the
study of continued fraction expansion, and also in some domains of
physics (crystallography) and computer science (formal
language theory, algorithms on words, pattern recognition).
Sturmian words admit several equivalent definitions and have numerous 
characterizations; in particular, they can be
characterized by their palindrome or return word structure
\cite{xDgP99pali, jJlV00retu}. For a comprehensive introduction to Sturmian words, see for instance \cite{jAjS03auto, jBpS02stur, nP02subs} and references therein. 

Sturmian words have exactly two factors of length 1,
and thus are infinite sequences over a two-letter alphabet $\cA = \{a,b\}$, say.
Here, an \emph{infinite word} (or \emph{sequence}) $\bx$ over $\cA$ is a map $\bx : \NN
\rightarrow \cA$. For any $i \geq 0$, we set $x_i = \bx(i)$ and
write $\bx = x_0x_1x_2\cdots$, each $x_i \in \cA$. Central
to our study is the following characterization of Sturmian words, which was
originally proved by Morse and Hedlund \cite{gHmM40symb}. An
infinite word $\bs$ over $\cA = \{a,b\}$ is Sturmian if and
only if there exists an irrational $\alpha \in (0,1)$, and a real
number $\rho$, such that $\bs$ is equal to one of the following two
infinite words:
\[
  s_{\alpha,\rho}, ~s_{\alpha,\rho}^{\prime}: \NN \rightarrow \cA
\]
defined by
\[
 \begin{matrix}
  &s_{\alpha,\rho}(n) = \begin{cases}
                        a    &\mbox{if} ~\lfloor(n+1)\alpha + \rho\rfloor -
                        \lfloor n\alpha + \rho\rfloor = 0, \\
                        b    &\mbox{otherwise};
                       \end{cases} \\
  &\qquad \\
  &s_{\alpha,\rho}^\prime(n) = \begin{cases}
                        a    &\mbox{if} ~\lceil(n+1)\alpha + \rho\rceil -
                        \lceil n\alpha + \rho\rceil = 0, \\
                        b    &\mbox{otherwise}.
                       \end{cases}
 \end{matrix} \qquad (n \geq 0)
\]
The irrational $\alpha$ is called the \emph{slope} of $\bs$ and
$\rho$ is the \emph{intercept}. If $\rho = 0$, we have
\[
  s_{\alpha,0} = ac_{\alpha} \quad \mbox{and} \quad s_{\alpha,0}^\prime =
  bc_{\alpha},
\]
where $c_{\alpha}$ is called the \emph{characteristic Sturmian
word} of slope $\alpha$ (see \cite{jBpS02stur}).

Our focus will be on palindromic factors of $c_\alpha$. In
general terms, a \emph{palindrome} is a finite word that reads the
same backwards as forwards. Palindromes are important tools used
in the study of factors of Sturmian words (e.g., \cite{aD81acom, aDfM94some, xD95pali, xDgP99pali}),  
and they have also become objects of great interest in computer science. 
The aim of this current paper is to completely describe where palindromes  
occur in $c_\alpha$ (and hence $s_{\alpha,0}$, $s_{\alpha,0}^\prime$). 
In order to do this, we shall make use of some previous results concerning factorizations of 
$c_\alpha$ into \emph{singular words}, which are particular palindromes. Singular words were first defined for the Fibonacci word $\bff$ (a special example of a Sturmian word) by Wen and Wen \cite{zWzW94some}, who established a decomposition of $\bff$ with respect to such words. This result was later extended by Melan\c{c}on \cite{gM99lynd}  to characteristic Sturmian words. More recently, Lev\'{e} and S\ee bold \cite{fLpS03conj} have generalized Wen and Wen's `singular' decomposition of $\bff$, by 
establishing a similar decomposition for each \emph{conjugate} of $\bff$ into what they called  
\emph{generalized singular words}. This last result has now been further extended by the present author   \cite{aG04conj} to $c_\alpha$ (and $c_{1-\alpha}$), where $\alpha$ has continued fraction expansion $[0;2,r,r,r,\ldots]$ for some $r\geq1$. 

It is well-known that any Sturmian word $\bs$ is 
\emph{uniformly recurrent}, i.e., any factor of $\bs$ occurs infinitely often
in $\bs$ with bounded gaps \cite{eCgH73sequ}. 
Accordingly, any palindromic factor $u$ of $c_\alpha$ has infinitely many 
occurrences in $c_\alpha$ and, as we shall see later (Corollary \ref{Cor:25.07.03(2)}), the distance between any two adjacent occurrences of $u$ is bounded above by an integer depending on $u$.   
Given any palindromic factor $u$ of $c_\alpha$, 
we shall 
establish a decomposition of $c_\alpha$ with respect to the occurrences of $u$. Such a decomposition shows precisely at which positions $u$ occurs in $c_\alpha$, and this is directly related to the continued fraction expansion of the irrational slope $\alpha$.

This paper is organized as follows. In Section \ref{S:prelim},
after some preliminaries on words and morphisms, we will recall 
some facts about $c_\alpha$ and consider some of its singular decompositions (Section \ref{SS:Sturmian}). Then, in Section 3, we consider the structure of palindromic factors of $c_\alpha$ with respect to its singular factors. We also recall the important notion of a return word and the concept of overlapping occurrences of a word in 
$c_\alpha$. Section 4 contains the lemmas we need  in order 
to establish the main result of this paper, which appears 
in Section 5.
Lastly, using results of Section 4, we obtain decompositions
of $c_\alpha$ that show precisely where a given factor of length
$q_n$ occurs in $c_\alpha$ (where $q_n$ is the denominator of the
$n$-th convergent to $\alpha = [0;1+d_1,d_2,d_3,\ldots]$, $d_i \geq 1$). 

\section{Preliminaries} \label{S:prelim}

Any of the following terminology that is not
further clarified can be found in either \cite{mL83comb} or
\cite{jBpS02stur}, which give more detailed presentations.

\subsection{Words and morphisms} \label{SS:words}

In what follows, let $\cA$ denote the two-letter alphabet $\{a,b\}$. A
(finite) \emph{word} is an element of the free monoid $\cAstar$
generated by $\cA$, in the sense of concatenation. The identity
$\empt$ of $\cAstar$ is called the \emph{empty word}, and
the \emph{free semigroup} over $\cA$
is defined by $\cAplus := \cAstar\setminus\{\empt\}$. We denote by
$\cAw$ the set of all infinite words over $\cA$, and define
$\cAinf := \cAstar \cup \cAw$. The \emph{length} $|w|$ of
a finite word $w$ is defined to be the number of letters it contains. 
(Note that $|\empt|=0$.)

A finite word $z$ is a \emph{factor} of a word $w \in \cAinf$ if $w =
uzv$ for some $u \in \cAstar$ and $v \in \cAinf$.
Furthermore, $z$ is called a \emph{prefix} (resp.~\emph{suffix})
of $w$ if $u = \empt$ (resp.~$v = \empt$), and we write $z \pref
w$ (resp.~$z \suff w$).  The word $z$ is said to have an \emph{occurrence} (or occur) at position 
$|u|$ of $w=uzv$, i.e., $z$ begins at the $|u|$-th position of $w$.
We denote by $|w|_{z}$ the number of occurrences of $z$ in $w$, i.e., the number of 
distinct positions at which $z$ occurs in $w$. For example,  
$|ababa|_{aba} = 2$ since $aba$ has two occurrences at positions 0 and 2 in $ababa$.

For any word $w \in \cAinf$, $\Omega(w)$ denotes the set of all
factors of $w$. Moreover, we denote by $\Omega_n(w)$ the set
of all factors of $w$ of length $n \in \NN$ (where $n \leq |w|$ 
for $w$ finite), i.e., $\Omega_n(w) = \Omega(w) \cap \cA^n$. If $u \in \Omega(w)$, then 
we shall simply write $u \prec w$.

The \emph{reversal operation} ~$\overset{\thicksim}{~~}$ ~in
$\cAstar$ is defined inductively by: $\rev{\empt} = \empt$ and,
for any $u \in \cAstar$ and $x \in \cA$, $(\rev{ux}) = x\rev{u}$.
Thus, if $w = x_0x_{1}x_{2}\ldots x_{n}$, with each $x_{i} \in \cA$,
then $\rev{w} = x_{n} x_{n-1}\ldots x_{1}x_0$. If $w = \rev{w}$, then
$w$ is called a \emph{palindrome}, and we define PAL to be the set
of all palindromes over $\cA$. It is useful to note that if $|w|$
is even, then $w$ is a palindrome if and only if $w = v\rev{v}$
for some word $v$. Otherwise, $w$ is a palindrome if and only if
$w = vx\rev{v}$ for some word $v$ and some letter $x \in \cA$.

The free monoid $\cAstar$ can be naturally embedded within a \emph{free
group}. We shall denote by $\cF$ the free group generated by $\cA$, which contains the 
\emph{inverse} $u^{-1}$ of each word $u \in \cAstar$. 
For any $u$, $v \in \cF$, we have $uu^{-1} = u^{-1}u = \empt$ and 
$(uv)^{-1} = v^{-1}u^{-1}$. If $u$, $w \in \cAstar$, we shall write 
$u^{-1}w$ (resp.~$wu^{-1}$) only if $u$ is a prefix (resp.~suffix) of $w$, so that 
$u^{-1}w$ (resp.~$wu^{-1}$) is a word in $\cAstar$. In particular, if $w=uv \in \cAstar$, then 
$u^{-1}w = v$ and $wv^{-1} = u$, and we have $|u^{-1}w| = |w|-|u|=|v|$, $|wv^{-1}| = |w|-|v| = |u|$. 

An \emph{endomorphism} (or simply \emph{morphism}) of  $\cAstar$ is a map $\psi: \cAstar \rightarrow
\cAstar$ such that $\psi(uv) = \psi(u)\psi(v)$ for all $u, v \in
\cAstar$.  It is uniquely determined by its image on the alphabet
$\cA$. Any morphism $\psi$ of $\cAstar$ can be uniquely extended to an
endomorphism of $\cF$ by defining $\psi(a^{-1}) = (\psi(a))^{-1}$
and $\psi(b^{-1}) = (\psi(b))^{-1}$, from which it follows that
$\psi(w^{-1}) = (\psi(w))^{-1}$ for any $w \in \cF$.

\subsubsection{Standard morphisms}

Define the following two morphisms of $\cAstar$:
\[
  E: \begin{array}{lll}
      a &\mapsto &b \\
      b &\mapsto &a
     \end{array}, \qquad \varphi: \begin{array}{lll}
                                   a &\mapsto &ab \\
                                   b &\mapsto &a
                                   \end{array}.
\]
A morphism $\psi$ of $\cAstar$ is \emph{standard} if $\psi(\bx)$ is a
characteristic Sturmian word for any characteristic Sturmian word
$\bx$ \cite{jBpS02stur}. In fact, a morphism $\psi$ is standard if and only if $\psi
\in \{E,\varphi\}^{*}$, i.e., if and only if it is a composition of
$E$ and $\varphi$ in any number and order \cite{aD97stan, jBpS02stur}.  
The standard morphisms $E$ and $\varphi E$ will play an important role in the proof of our main result.

\subsection{Characteristic Sturmian words $c_\alpha$ and singular words}
\label{SS:Sturmian}

Note that every irrational $\alpha \in (0,1)$ has a unique
continued fraction expansion
\[
  \alpha = [0;a_1,a_2,a_3,\ldots] = \cfrac{1}{a_1+
                                \cfrac{1}{a_2 +
                                \cfrac{1}{a_3 + \cdots
                                 }}}
\]
where each $a_i$ is a positive integer. If the sequence
$(a_i)_{i\geq1}$ is eventually periodic, with $a_i = a_{i+m}$ for
all $i \geq n$, we use the notation $\alpha =
[0;a_1,a_2,\ldots,a_{n-1},\ov{a_{n},a_{n+1},\ldots,a_{n+m-1}}].$
The \emph{$n$-{th} convergent} to $\alpha$ is defined by
\[
  \frac{p_n}{q_n} = [0;a_1,a_2,\ldots,a_n], \quad \mbox{for all}~ n\geq 1,
\]
where the sequences $(p_n)_{n\geq0}$ and $(q_n)_{n\geq0}$ are
given by
\[
\begin{matrix}
 &p_{0} = 0, &p_{1} = 1, &p_n = a_np_{n-1} + p_{n-2},  &n\geq 2; \\
 &q_{0} = 1, &q_{1} = a_1, &q_n = a_nq_{n-1} + q_{n-2}, &n\geq 2.
\end{matrix}
\]

Suppose $\alpha = [0;1+d_1,d_2,d_3, \ldots]$ with $d_1 \geq 0$
and all other $d_n > 0$. To the \emph{directive sequence}
$(d_1,d_2,d_3,\ldots)$, we associate a sequence $(s_n)_{n \geq
-1}$ of words defined by
\[
  s_{-1} = b, ~s_{0} = a, ~s_{n} = s_{n-1}^{d_{n}}s_{n-2}, \quad n \geq 1.
\]
Such a sequence of words is called a \emph{standard sequence}, and
we have
\[
  |s_n| = q_n \quad \mbox{for all}~ n\geq0.
\]
Note that $ab$ is a suffix of $s_{2n-1}$ and $ba$ is a suffix of
$s_{2n}$, for all $n \geq 1$.

Standard sequences are related to characteristic Sturmian words in
the following way. Observe that, for any $n\geq0$, $s_n$ is a
prefix of $s_{n+1}$, which gives obvious meaning to ${\lim}_{n
\rightarrow \infty}s_n$ as an infinite word. In
fact, each $s_n$
is a prefix of $c_\alpha$, and we have 
\begin{equation} \label{eq:ACW1}
c_{\alpha} = \underset{n \rightarrow \infty}{\mbox{lim}}s_n \quad \mbox{(see \cite{aFmMuT78dete, tB93desc})}.  
\end{equation}

\subsubsection{Some singular decompositions of $c_\alpha$}

Note that if $\alpha = [0;1,d_1,d_2,d_3, \ldots]$,
then
\begin{equation} \label{eq:observe}
  1 - \alpha = \frac{1}{1+1/(1/\alpha-1)} = [0;1+d_1,d_2,d_3,\ldots].
\end{equation}
For any irrational $\alpha \in (0,1)$, $E(c_\alpha) =
c_{1-\alpha}$, i.e., $c_{1-\alpha}$ is obtained from $c_\alpha$ by
exchanging $a$'s and $b$'s \cite{bP97prop}. Thus, in light of the
above observation \eqref{eq:observe}, we shall hereafter restrict our 
attention to the case when $\alpha = [0;1+d_1,d_2,d_3,\ldots]$ with $d_1\geq1$.

Melan\c{c}on \cite{gM99lynd} (also see \cite{wCzW03some, zWzW94some}) has 
introduced the singular words $(w_n)_{n\geq0}$ of $c_\alpha$  defined by
\[
  w_n = \begin{cases}
         as_{n}b^{-1} &\mbox{if $n$ is odd}, \\
         bs_na^{-1} &\mbox{otherwise}.
        \end{cases}
\] 
Moreover, for each $n \geq -1$, Melan\c{c}on \cite{gM99lynd} defined the words 
\[
  v_n = \begin{cases} as_{n+1}^{d_{n+2} - 1}
                          s_n b^{-1} &\mbox{if $n$ is odd}, \\
                          bs_{n+1}^{d_{n+2} - 1}
                          s_n a^{-1} &\mbox{otherwise}.
                          \end{cases}
\]
Clearly, the word $v_n$ differs from $w_{n+2}$ by a factor
$s_{n+1}$, and it is easily proved that all $v_n$ and $w_n$ are
palindromes. Here, we will call $w_n$ (resp.~$v_n$) the \emph{$n$-th singular word} (resp.~\emph{$n$-{th} adjoining singular word}) of $c_{\alpha}$,
and use the convention $w_{-2}=v_{-2} = \empt$, $w_{-1} = a$.

Singular words play an important role in the study of factors of
Sturmian words. In particular, as we shall see in the next
section, the words $w_n$ and $v_{n-1}$ can be used to determine
the structure of all palindromic factors of a Sturmian word of
slope $\alpha$. We have the following decomposition of
$c_\alpha$ in terms of singular and adjoining singular words.

\begin{proposition} \emph{\cite{zWzW94some, gM99lynd}} \label{P:melancon1}
~$c_\alpha = \prod_{j=-1}^\infty(v_{2j}w_{2j+1})^{d_{2j+3}}=
        \prod_{j=-1}^\infty v_{j}$.\qed
\end{proposition}

\begin{notation} In order to simplify proceedings, we introduce some notation.
\begin{itemize}
\item[(i)]
Let $\gamma \in (0,1)$ be irrational with $\gamma =
[0;a_1,a_2,a_3,\ldots]$.  For any $n \in \NN$ and 
integer $k$ such that $k \geq 1-a_{n+1}$, define
\[
  \gamma_{n,k} := [0;a_{n+1}+k,a_{n+2},a_{n+3},\ldots] 
 \]
and write $\gamma_{n,0} = \gamma_n$. Note that 
$\gamma_{0,0} = \gamma = [0; a_1,a_2,\ldots,a_n + \gamma_{n}]$ for all $n\geq1$.
\item[(ii)] As $c_\alpha$ is uniformly recurrent, given any factor
$w$ of $c_\alpha$, the occurrences of $w$ in $c_\alpha$ can be 
arranged as a sequence $(w^{(i)})_{i\geq 1}$, where $w^{(i)}$ denotes the
$i$-{th} occurrence of $w$ in $c_\alpha$.
\end{itemize}
\end{notation}

With the above notation, we may now state a corollary of Proposition 
\ref{P:melancon1}.
\newpage
\begin{corollary} \label{Cor:Melancon2} 
Let $n \in \NN$ be fixed. The characteristic Sturmian word
$c_\alpha$ has the following two decompositions:
\begin{itemize}
\item[$(1)$] $\qquad$ \vspace{-0.4cm}
\[
  c_\alpha = \left(\prod_{j=-1}^{n-1}(v_{2j}w_{2j+1})^{d_{2j+3}}\right)
      w_{2n}^{(1)}z_{1}
      w_{2n}^{(2)}z_2w_{2n}^{(3)}z_3
      \cdots,
\] 
where $\bz := z_1z_2z_3\cdots$ is given by $c_{\alpha_{2n+1}}$ over the alphabet $\{v_{2n-1}, w_{2n+1}\}$.
\item[$(2)$] $\qquad$ \vspace{-0.4cm}
\[
  c_\alpha = \left(\prod_{j=-1}^{n-1}(v_{2j}w_{2j+1})^{d_{2j+3}}\right)z_{1}
      W_{1}z_2W_2z_3
      \cdots,
\] where $\bz := z_1z_2z_3\cdots$ is given by $c_{\alpha_{2n+1,1}}$ over the alphabet 
$\{w_{2n}, v_{2n-2}\}$ and, for all $i\geq 1$, 
\[
  W_i   = \begin{cases}
                                      v_{2n-1} &\mbox{if} ~z_{i}=w_{2n}, \\
                                      w_{2n-1} &\mbox{if} ~z_{i}=v_{2n-2}.
                              \end{cases}
\]
\end{itemize} 
\end{corollary}
\begin{proof}
See \cite[Corollary 4.6]{gM99lynd}.
\end{proof}

\begin{example} \label{Ex:Fibonacci}
The best known example of a characteristic Sturmian word is the
infinite \emph{Fibonacci word} $\bff$, which has been extensively
studied by many authors (see \cite{aD81acom, xD95pali}, for
example). It is well-known that
\[
 \bff = \underset{n \rightarrow \infty}{\mbox{lim}} f_n =
abaababaabaababaababaabaababaabaab \cdots,
\]
where $(f_n)_{n\geq-1}$ is the sequence of \emph{finite Fibonacci
words} defined by
\[
  f_{-1} = b, ~f_{0} = a, ~f_n = f_{n-1}f_{n-2}, \quad n \geq 1.
\]
Clearly, $|f_{n}| = F_n$, where $F_n$ is the \emph{$n$-{th}
Fibonacci number} defined by
\[
  F_{-1} = 1, ~F_{0}=1, ~F_{n} = F_{n-1} + F_{n-2}, \quad n \geq 1.
\]
Note that $(f_n)_{n\geq-1}$ is a standard sequence with associated
directive sequence $(1,1,1,\ldots)$, and hence $w_n = v_n$ for all $n
\geq -1$. Moreover, in view of
\eqref{eq:ACW1}, $\bff = c_\alpha$ where 
$\alpha = (3 - \sqrt{5})/2 = [0;2,\overline{1}]$, in which case 
$\alpha = \alpha_{2n+1,1}$ and $1 - \alpha = \alpha_{2n+1}$, for all $n \in \NN$.
Hence, $\bff = c_{\alpha_{2n+1,1}}$ and 
$E(\bff) = c_{\alpha_{2n+1}}$. Accordingly, one deduces from the  
above corollary that
\[
  \bff = \left(\prod_{j=-1}^{n-1}w_j\right)w_n^{(1)}z_1w_n^{(2)}z_2
w_n^{(3)}z_3\cdots,
\]
where $\bz := z_1z_2z_3\cdots$ is the
Fibonacci word over the alphabet $\{w_{n+1}, w_{n-1}\}$  
(also see \cite[Theorem 2]{zWzW94some}). For instance, when $n = 2$, $w_{n-1} = w_1 = aa$, $w_{n} = w_2 = bab$, $w_{n+1} = w_3 = aabaa$,
and $\bz$ is the Fibonacci word over the
alphabet $\{aabaa, aa\}$. Indeed, one may write
\[
  \bff = abaa(bab)aabaa(bab)aa(bab)aabaa(bab)aabaa(bab)aa(bab)aabaa(bab)aa(bab)
      aabaa\cdots.
\] 
\end{example}

\section{Palindromes, return words and overlap}

\subsection{Structure of palindromes in $c_\alpha$} \label{S:palindromes}

In \cite{wCzW03some}, Cao and Wen considered the structure of palindromic factors of 
$c_\alpha$ with respect to singular words. Specifically, they proved the following result 
concerning palindromic factors $u$ of $c_\alpha$ with $q_n < |u| \leq q_{n+1}$.
(For technical reasons, we set $q_{-1} = 1$, so that $|s_n| = q_n$ for all $n
\geq -1$.)
\newpage
\begin{proposition} \label{P:11} \emph{\cite{wCzW03some}}
Let $u \in$ \emph{PAL} with $q_n < |u| \leq q_{n+1}$ for some $n
\in \NN$. Then $u \prec c_\alpha$ if and only if $u$ takes one of
the following forms:
\begin{itemize}
\item[$(1)$] $u = vw_n\rev{v}$ with $v \suff v_{n-1}$ and $|v| \leq
\frac{1}{2}|v_{n-1}| = \frac{1}{2}(q_{n+1}-q_n);$
\item[$(2)$] $u = vv_{n-1}\rev{v}$ with $v \suff w_{n}$ and $|v| \leq
\frac{1}{2}q_n;$
\item[$(3)$] $u = v(w_{n-1}v_{n-2})^kw_{n-1}\rev{v}$ with $v \suff
v_{n-2}$, $v \ne v_{n-2}$, and $0\leq k\leq d_{n+1} - 2;$ 
\item[$(4)$] $u = v(v_{n-2}w_{n-1})^kv_{n-2}\rev{v}$ with $v \suff
w_{n-1}$, $v \ne w_{n-1}$, and $0\leq k\leq d_{n+1} - 1;$ 
\item[$(5)$] $u = w_{n+1}$. 
\end{itemize}
Moreover, if $k=0$ in $(3)$ $($resp.~$(4))$, then $|v| > \frac{1}{2} |v_{n-2}| = \frac{1}{2}(q_n - q_{n-1})$ 
$($resp.~$|v| > \frac{1}{2} q_{n-1})$. \qed
\end{proposition}

Hereafter, we will make frequent use of  
the following properties of singular words. Some of these properties may be used without referring to the given lemma.

\begin{lemma} \label{L:2} \emph{\cite{gM99lynd, wCzW03some}}
Let $x$, $y \in \cA$ $(x \ne y)$ with $y \suff s_n$. Then, for any
$n \in \NN$,
\begin{itemize}
\item[$(1)$] ~$yx^{-1}w_n = ys_ny^{-1} = w_{n-1}v_{n-2}, \quad 
w_nx^{-1}y = v_{n-2}w_{n-1};$
\item[$(2)$] ~$w_{n+1} = w_{n-1}v_{n-2}v_{n-1}= v_{n-1}v_{n-2}w_{n-1};$
\item[$(3)$] ~$v_{n-1} = (w_{n-1}v_{n-2})^{d_{n+1}-1}w_{n-1};$
\item[$(4)$] ~$w_{n+1} = (w_{n-1}v_{n-2})^{d_{n+1}}w_{n-1};$
\item[$(5)$] ~$w_{n+1} = y\prod_{j=-1}^{n-1}v_j;$
\item[$(6)$] ~$w_{n} \nprec w_{n+1};$
\item[$(7)$] ~$v_{n-1} \nprec w_n$.
\end{itemize} \qed
\end{lemma}

Now, for each $n\in\NN$ and $0\leq k \leq d_{n+1}-1$, let us denote by $U_{n,k}$ and 
$\ov{U}_{n,k}$ the palindromes given by
\[
  U_{n,k}  := (w_{n-1}v_{n-2})^kw_{n-1} \quad  \mbox{and} \quad \ov{U}_{n,k} := (v_{n-2}w_{n-1})^kv_{n-2}.  
\]  
Note that $U_{n,k} = w_{n-1}\ov{U}_{n,k}(v_{n-2})^{-1}$. 
Also observe that  the singular words $(w_{n})_{n\geq-1}$ and $(v_{n})_{n\geq-1}$ are given by
\[
  w_{n-1} = U_{n,0}  \quad \mbox{and} \quad v_{n-1} = U_{n,d_{n+1}-1} = \ov{U}_{n+1,0} \quad 
  \mbox{for all $n\geq0$}.
\]  

From the preceding proposition and Lemma \ref{L:2}, we easily deduce the following result, which gives the structure of all palindromic factors of $c_\alpha$ in terms of $U_{n,k}$ and $\ov{U}_{n,k}$. 
The proof is left to the reader.

\begin{corollary} \label{Cor:all_palindromes}
Let $u \in$ \emph{PAL} with $|u|\geq2$. Then $u$ is a factor of $c_\alpha$ if and only if, for some $n\in\NN$, we have 
\begin{equation} \label{eq:U1}
\hspace{-4cm}
  u = vU_{n,k}\rev{v}, ~\mbox{where $v \suff 
                                                                 v_{n-2}$, $v\ne v_{n-2}$ and $0\leq k \leq d_{n+1}-2$} 
\end{equation}
or 
\begin{equation} \label{eq:U2}
\hspace{-3.9cm}
  u = v\ov{U}_{n,k}\rev{v}, ~\mbox{where $v \suff 
                                                                  w_{n-1}$, $v\ne w_{n-1}$ and $0\leq k \leq d_{n+1}-1$}.
\end{equation} \qed
\end{corollary}
\begin{note} Let us point out that $\ov{U}_{0,0} = \empt$ and $U_{0,k-1} = a^k = \ov{U}_{0,k}$ for $1\leq k \leq d_{1}-1$. Therefore, if $u$ takes the form \eqref{eq:U1} or \eqref{eq:U2} for $n=0$, then $u = a^k$ 
for some $k \in [2,d_1-1]$. So, if $d_1\leq 2$, a palindromic factor of $c_\alpha$ is given by \eqref{eq:U1}  or \eqref{eq:U2} for some $n\geq1$.
\end{note}

\begin{remark} It is important to note that Corollary \ref{Cor:all_palindromes} (and 
also Proposition \ref{P:11}) gives
the structure of all palindromic factors of any Sturmian word of
slope $\alpha$. Indeed, Mignosi \cite{fM89infi} proved that any
two Sturmian words $\bs$, $\bt$ of the same slope are
\emph{equivalent}, i.e., $\Omega(\bs) = \Omega(\bt)$. Whence, for
any real number $\rho$, we have
\[
  \Omega(s_{\alpha,\rho}) = \Omega(s_{\alpha,\rho}^\prime) 
  = \Omega(c_\alpha),
\]
i.e., a palindrome is a factor of some Sturmian word of slope $\alpha$ if and only if it is a factor of 
$c_\alpha$.
\end{remark}

\subsection{Return words and overlapping occurrences} \label{S:return}

Let us write $c_\alpha = x_0x_1x_2\cdots$, each $x_i \in \cA$, and
let $w \prec c_\alpha$. Suppose $n_1 < n_2 < n_3 < \cdots$ are all
the natural numbers $n_i$ such that $w = x_{n_i}x_{n_i+1}\cdots
x_{n_i+|w|-1}$. Then the word $x_{n_i}\cdots x_{n_{i+1}-1}$ is a
\emph{return word} of $w$ in $c_\alpha$. That is, we define the set $\cR_w(c_\alpha)$ of
return words of $w$ to be the set of all distinct words beginning
with an occurrence of $w$ and ending exactly before the next
occurrence of $w$ in $c_\alpha$. This notion was introduced
independently by Durand \cite{fD98acha}, and Holton and Zamboni
\cite{cHlZ99desc}.  Clearly,
$\cR_w(c_\alpha)$ is finite since 
the distance between two adjacent occurrences of $w$ in
$c_\alpha$ is bounded. In fact, Vuillon \cite{lV01acha} has proved that an
infinite word $\bs$ over $\cA$ is Sturmian if and only if, for any
factor $w$ of $\bs$, there are exactly two return words of $w$ in
$\bs$. Suppose $\cR_w(c_\alpha) = \{u_1,u_2\}$. Then $c_\alpha$ can be
uniquely factorized as $c_\alpha = vu_{i_1} u_{i_2}\cdots
u_{i_k}\cdots$, where each $i_k \in \{1,2\}$ and the first occurrence of $w$ in 
$c_\alpha$ is at position $|v|$. The
infinite word $\cD_w(c_\alpha) := u_{i_1}u_{i_2}\cdots u_{i_k}\cdots$,
called the \emph{derived word of $c_\alpha$ with respect to $w$},
can be viewed as an infinite word over the alphabet $\{u_1,u_2\}$. In particular, $\cD_w(c_\alpha)$ is a Sturmian word over
the alphabet $\cR_w(c_\alpha)$ \cite{iFlV02gene}.
For example, the return words of
$w_n$ in $\bff$ are $w_nw_{n+1}$ and $w_{n}w_{n-1}$, and $\cD_{w_n}(\bff)$ is the 
Fibonacci word over the alphabet $\{w_nw_{n+1},w_{n}w_{n-1}\}$ (see Example \ref{Ex:Fibonacci}).

Given $w \prec c_\alpha$, a return word of $w$ in $c_\alpha$ is
not necessarily longer than $w$, in which case $w$ has overlapping
occurrences in $c_\alpha$.  More precisely, if there exist
non-empty words $u$, $v$ and $z$ such that $w = uz = zv$ and $uzv
\prec c_\alpha$, then $w$ is said to have {\em overlap} in
$c_\alpha$ with \emph{overlap factor} $z$. Further, one can write 
$uzv = wz^{-1}w$; whence $w$ has overlap in $c_\alpha$ if $wz^{-1}w \prec
c_\alpha$ for some $z \in \cA^+$. In this case, $wz^{-1}$ is a
return word of $c_\alpha$ that has length less than that of $w$. Clearly,
since any factor $w$ of $c_\alpha$ has exactly two return words, $w$ 
has at most two different overlap factors.

Return words, and the concept of overlap, are fundamentally important to our 
study of
occurrences of palindromes in $c_\alpha$. Indeed, we shall be
establishing decompositions of $c_\alpha$ with respect to certain palindromic 
factors that have overlap, i.e., 
palindromic factors $u$ that have a return word (or return words) of
length(s) less than $|u|$. Specifically, given any palindromic factor 
$u$ of $c_\alpha$, we can write
\[
  c_\alpha = z_0u^{(1)}z_1u^{(2)}z_2u^{(3)}z_3\cdots,
\]
where $z_0 \in \cAstar$ and all other $z_i$ are such that $z_i^{-1} \in \cA^+$ or $z_i \in \cAstar$,  according to whether the occurrences $u^{(i)}$ and $u^{(i+1)}$ do or do not overlap each other,  respectively. For instance, if $u = w_n$ is the $n$-th singular factor of the Fibonacci word, then, as shown in  Example \ref{Ex:Fibonacci}, each $z_i \in \{w_{n+1}, w_{n-1}\}$ ($i\geq 1$); in which case $u$ does not have overlap in $\bff$. 

The following result shows precisely which factors of $c_\alpha$
have no overlapping occurrences in $c_\alpha$.

\begin{proposition} \emph{\cite[Theorem 10]{wCzW03some}} \label{P:bCzW03-10}
Let $u \prec c_\alpha$ with $q_n < |u| \leq q_{n+1}$ for some $n
\in \NN$. Then $u$ has no overlap in $c_\alpha$ if and only if $u
= w_{n+1}$, or $w_n \prec u$. \qed
\end{proposition}

Accordingly, one easily deduces from Proposition \ref{P:11} and Lemma \ref{L:2} which palindromic factors of $c_\alpha$ do not have overlap. 

\begin{corollary} \label{Cor:7}
Let $u \in$ \emph{PAL} and $u \prec c_\alpha$ with $q_n < |u| \leq q_{n+1}$ 
for some $n \in \NN$. Then $u$ is a palindrome without overlap in 
$c_\alpha$ if and only if $u = w_{n+1}$, or $u = vw_n\rev{v}$ with 
$v \suff v_{n-1}$ and $|v| \leq \frac{1}{2}|v_{n-1}|$. \qed
\end{corollary}

\section{Decompositions of $c_\alpha$ into palindromes}

In this section, we prove some lemmas which lead us to the main result 
of this paper (Theorem \ref{T:21.08.03(1)}). 

\subsection{Useful results} \label{SS:useful_results}

In what follows, let us denote by $G$ the standard morphism of $\cAstar$ given by
\[
  G = \varphi E :  \begin{array}{lll}
                               a &\mapsto &a \\
                               b &\mapsto &ab
                              \end{array}.
\]

\begin{lemma} \label{L:parvaix97} \emph{\cite{bP97prop}}
For any irrational $\gamma \in (0,1)$,
  $E(c_\gamma) = c_{1 - \gamma}$ and $G(c_\gamma) = c_{\gamma/
(1+\gamma)}$. \qed
\end{lemma}

The following simple, yet useful, corollary (and the remark to follow) will be needed in our proofs.

\begin{corollary} \label{Cor:14.08.03(2)}
For any irrational $\gamma \in (0,1)$ and $k \in \NN$,
$G^k(c_\gamma) = c_{\gamma/(1+k\gamma)}$.
\end{corollary}
\begin{proof} Induction on $k$.
\end{proof}
\begin{remark} \label{R:G^k_and_alpha} Recall that we are restricting our attention to $c_\alpha$  
where $\alpha$ has continued fraction expansion $[0;1+d_1,d_2,d_3,\ldots]$, 
$d_1\geq 1$. Let us note that 
$\frac{\alpha}{1 + k\alpha} = \frac{1}{k + 1/\alpha} =
[0;1+d_1+k,d_2,d_3,\ldots] = \alpha_{0,k}$ and, more generally, $\frac{\alpha_n}{1 + k\alpha_n} = 
 [0;d_{n+1}+k,d_{n+2},d_{n+3},\ldots]$ for $n\geq 1$. Consequently, 
$$G^{k}(c_{\alpha_{n}}) = c_{\alpha_{n,k}} \quad \mbox{for all $n\geq0$}.$$
It is also easily checked that $1 - \alpha_{n+1,1} =
[0;1,d_{n+2},d_{n+3},\ldots] = \alpha_{n,1-d_{n+1}}$, for any $n\geq1$; whence 
\begin{equation} \label{eq:exchange}
E(c_{\alpha_{n+1,1}}) = c_{\alpha_{n,1-d_{n+1}}}  \quad \mbox{for all $n \geq 1$}.
\end{equation}
(Note that $E(c_{\alpha_{1,1}}) = c_{\alpha_{0,-d_1}}$.)
\end{remark}

\subsection{Some lemmas}

Here, we simplify Melan\c{c}on's decompositions of $c_\alpha$,  
given in Corollary \ref{Cor:Melancon2}. In particular, we obtain two different decompositions of
$c_\alpha$ with respect to occurrences of the palindromes 
\[
  U_{n,k}  = (w_{n-1}v_{n-2})^kw_{n-1} \quad  \mbox{and} \quad \ov{U}_{n,k} = (v_{n-2}w_{n-1})^kv_{n-2}  \quad (0 \leq k \leq d_{n+1}-1), 
\]  
which form the basis of all palindromic factors of $c_\alpha$ (see Corollary 
\ref{Cor:all_palindromes}). 
From the first of these decompositions,  we easily deduce decompositions of $c_\alpha$ that show exactly where the singular words $w_n$ and $v_n$ occur in $c_\alpha$, for any $n \in \NN$. (Recall that Corollary \ref{Cor:Melancon2} gives a decomposition of $c_\alpha$ which shows all of the occurrences of $w_{2n}$, but this result does not provide information as to the exact positions of $w_{2n-1}$ in $c_\alpha$.) 

\begin{notation} For any morphism $\psi$ of $\cAstar$ such that $\psi(a) = u$ and $\psi(b) = v$ for some $u$, $v \in \cAstar$, we shall write $\psi = (u,v)$ to indicate the image of $\psi$ on the alphabet $\cA$ .  If $\bx= x_0x_1x_2\cdots \in \cAw$, then $\psi(\bx) =
\psi(x_0)\psi(x_1)\psi(x_2)\cdots$ is the word obtained from $\bx$ by
replacing the letters $a$ and $b$ in $\bx$ by the words $u$ and $v$,
respectively. We shall denote by $\bx\{u,v\}$ the word $\psi(\bx)$.  
In particular, $c_\alpha\{u,v\}$ denotes the characteristic
Sturmian word of slope $\alpha$ over the alphabet $\{u,v\}$.
\end{notation}

\begin{lemma} \emph{\cite{wCzW03some}} \label{L:21.08.03(2)}
For any $n \geq 1$, $c_\alpha = c_{\alpha_{n,1}}
\{s_n,s_{n-1}\}$. \qed
\end{lemma}

\begin{lemma} \label{L:14.08.03(4)} For any $n \in \NN$,
\[
\prod_{j=-1}^{n-1}(v_{2j}w_{2j+1})^{d_{2j+3}} = \prod_{j=-1}^{2n-1}v_{j}.
\]
\end{lemma}
\begin{proof} Using Lemma \ref{L:2}$(3)$, observe that
for any integer $j \geq 0$,
\[
(v_{2j}w_{2j+1})^{d_{2j+3}} =
v_{2j}w_{2j+1}(v_{2j}w_{2j+1})^{d_{2j+3}-1}
                            =  v_{2j}(w_{2j+1}v_{2j})^{d_{2j+3}-1}
                                w_{2j+1}
                            = v_{2j}v_{2j+1},
\]
from which the result is readily deduced.
\end{proof}

\begin{lemma} \label{L:cool2}
For any $n\geq1$ and $0 \leq k \leq d_{n+1} - 1$,  we have 
\[
  c_\alpha = \left(\prod_{j=-1}^{n-2}v_j\right)U_{n,k}^{(1)}z_1U_{n,k}^{(2)}z_2
U_{n,k}^{(3)}z_3\cdots,
\]
where $\bz := z_1z_2z_3\cdots$ is given by $c_{\alpha_{n,-k}}$ over the alphabet 
$\{(U_{n,k-1})^{-1}, w_{n}\}$.
\end{lemma}
\begin{note} We set $U_{n,-1} = (v_{n-2})^{-1}$ and $\ov{U}_{n,-1} = (w_{n-1})^{-1}$; whence if  $k = 0$, then 
\begin{equation*} 
  (U_{n,k-1})^{-1} = (U_{n,-1})^{-1} = v_{n-2}   \quad \mbox{and} \quad 
  (\ov{U}_{n,k-1})^{-1} = (\ov{U}_{n,-1})^{-1} = w_{n-1}. 
\end{equation*}
\end{note}
\begin{proof}[Proof of Lemma $\ref{L:cool2}$] We first prove the result  for odd $n = 2m+1$, $m\geq0$. By Corollary \ref{Cor:Melancon2}$(1)$ and Lemma \ref{L:14.08.03(4)}, we have
\[
 c_\alpha = \left(\prod_{j=-1}^{2m-1}v_j\right)\psi(c_{\alpha_{2m+1}}),
\]
where $\psi = (w_{2m}v_{2m-1}, w_{2m}w_{2m+1})$. Further, by Remark
\ref{R:G^k_and_alpha}, we have 
\[
  \psi G^{k}(c_{\alpha_{2m+1,-k}}) = \psi(c_{\alpha_{2m+1}}) \quad \mbox{for $0 \leq k \leq d_{2m+2}-1$}.
\]
Therefore, 
\[
 c_\alpha =
 \left(\prod_{j=-1}^{2m-1}v_{j}\right)\psi G^k(c_{\alpha_{2m+1,-k}}),
 \]
where $G^k = (a,a^kb)$, and hence
\begin{align*}
  \psi G^{k} 
         &= (w_{2m}v_{2m-1}, (w_{2m}v_{2m-1})^{k}w_{2m}
               w_{2m+1}) \\
         &= (w_{2m}v_{2m-1}, U_{2m+1,k}w_{2m+1} ) \\
         &= (U_{2m+1,k}(U_{2m+1,k-1})^{-1}, U_{2m+1,k}w_{2m+1} ).
\end{align*}
Clearly, $U_{2m+1,k}w_{2m+1}$ and
$w_{2m}v_{2m-1}$ ($=U_{2m+1,k}(U_{2m+1,k-1})^{-1}$) must be the two 
return words of $U_{2m+1,k}$. Also, using Lemma \ref{L:2}, we find that 
$U_{2m+1,k}$ is not a factor of the prefix $(\prod_{j=-1}^{2m-1}v_{j})U_{2m+1,k-1}$ of $c_\alpha$, since 
\[
  \left(\prod_{j=-1}^{2m-1}v_{j}\right)U_{2m+1,k-1} = x^{-1}w_{2m}v_{2m-1}U_{2m+1,k-1} = 
   x^{-1}U_{2m+1,k} \quad (x \in \cA). 
\]  
Thus, the derived word of $c_\alpha$ with respect to $U_{2m+1,k}$ is given 
by $\cD_{U_{2m+1,k}}(c_\alpha) = c_{\alpha_{2m+1,-k}}\{w_{2m}v_{2m-1}, U_{2m+1,k}w_{2m+1}\}$, and we can write 
\[
c_\alpha =
\left(\prod_{j=-1}^{2m-1}v_j\right)U_{2m+1,k}^{(1)}z_1U_{2m+1,k}^{(2)}
              z_2U_{2m+1,k}^{(3)}z_3\cdots,
\]
where $\bz := z_1z_2z_3\cdots$ is given by $c_{\alpha_{2m+1,-k}}$ over the alphabet 
$\{(U_{2m+1,k-1})^{-1}, w_{2m+1}\}$. This completes the proof for odd $n$.

Let us now prove that the assertion holds for even $n=2m$, $m\geq 1$. By considering occurrences of $U_{2m-1,d_{2m}-1} $ ($= v_{2m-2}$) in $c_\alpha$, one deduces from the above that, for any integer $m \geq 1$,
\begin{equation} \label{eq:020605}
  c_\alpha = \left(\prod_{j=-1}^{2m-3}v_j\right) v_{2m-2}\phi(c_{\alpha_{2m-1,1-d_{2m}}}),
\end{equation}
where 
$\phi = ((U_{2m-1,d_{2m}-2})^{-1}v_{2m-2}, w_{2m-1}v_{2m-2})$.              
In fact, using $(2)$ and $(3)$ of Lemma \ref{L:2}, we can write 
$\phi = (v_{2m-3}w_{2m-2}, w_{2m-1}v_{2m-2}) =
        ((v_{2m-2})^{-1}w_{2m}, w_{2m-1}v_{2m-2})$.
Again, using Remark \ref{R:G^k_and_alpha}, we have 
$EG^{k+1}(c_{\alpha_{2m,-k}}) = E(c_{\alpha_{2m,1}}) = c_{\alpha_{2m-1,1-d_{2m}}}$,
where $EG^{k+1} = (b,b^{k+1}a)$. Whence, it follows from \eqref{eq:020605} that 
\[
  c_\alpha =
  \left(\prod_{j=-1}^{2m-2}v_j\right)\phi E G^{k+1}(c_{\alpha_{2m,-k}}),
\]
where 
\[
  \phi E G^{k+1}  
  =   (w_{2m-1}v_{2m-2}, (w_{2m-1}v_{2m-2})^{k+1}(v_{2m-2})^{-1}w_{2m})  
  = (U_{2m,k}(U_{2m,k-1})^{-1}, U_{2m,k}w_{2m}).
\] 
The result now follows (as for the odd case) since $U_{2m,k}w_{2m}$ and $w_{2m-1}v_{2m-2}$ ($= U_{2m,k}(U_{2m,k-1})^{-1}$) are the two return words of $U_{2m,k}$. 

\end{proof}

\begin{remark}
From Lemma \ref{L:cool2},  we readily deduce two `singular' decompositions of $c_\alpha$ with respect to the occurrences of $w_{n}$ and $v_{n}$, for any $n\in\NN$. Indeed,  we have $U_{n,0} = w_{n-1}$ and 
$U_{n,d_{n+1}-1} = v_{n-1}$. Therefore, taking $k=0$ in the above lemma, we obtain a decomposition that shows exactly where the $n$-th singular word $w_n$ occurs in $c_\alpha$. That is, 
for any $n\geq0$, 
 \begin{equation} \label{eq:w_n}
  c_\alpha = \left(\prod_{j=-1}^{n-1}v_{j}\right)
       w_{n}^{(1)}z_{1}w_{n}^{(2)}z_2w_{n}^{(3)}z_3
      \cdots,
\end{equation}      
where $\bz := z_1z_2z_3\cdots$ is given by $c_{\alpha_{n+1}}$ over the alphabet 
$\{v_{n-1}, w_{n+1}\}$. 

Now, taking $k = d_{n+1}-1$, we find that, for any $n\geq 0$,    
\begin{equation} \label{eq:v_{n-1}}
  c_\alpha  = \left(\prod_{j=-1}^{n-1}v_j\right)v_{n}^{(1)}z_1v_{n}^{(2)}z_2
v_{n}^{(3)}z_3\cdots,
\end{equation}
where $\bz := z_1z_2z_3\cdots$ is given by $c_{\alpha_{n+2,1}}$ over the alphabet 
$\{w_{n+1}, (U_{n+1,d_{n+2}-2})^{-1}\}$ (since \\ $E(c_{\alpha_{n+1,1-d_{n+2}}}) = c_{\alpha_{n+2,1}}$).
This also holds for $n=-1$ since, from Lemma \ref{L:21.08.03(2)},  we have
\[
  c_{\alpha} = c_{\alpha_{1,1}}\{s_1,s_0\} = c_{\alpha_{1,1}}\{a^{d_1}b,a\} = c_{\alpha_{1,1}}\{v_{-1}w_0, w_{-1}v_{-2}\}.
\]
\end{remark}

The following simple decomposition of $c_\alpha$ (which has also been proved independently in 
\cite{wCzW03some}) is a direct consequence of \eqref{eq:v_{n-1}}. 

\begin{proposition} \label{P:21.08.03(3)}
For any $n \in \NN$, we have
\[
  c_\alpha = \left(\prod_{j=-1}^{n-2}v_j\right)c_{\alpha_{n+1,1}}\{v_{n-1}w_{n},w_{n-1}v_{n-2}\}.
\]
\qed
\end{proposition}

\begin{lemma} \label{L:cool3}
For any $n \geq 1$ and 
$0 \leq k \leq d_{n+1} - 1$, we have 
\[
  c_\alpha = \left(\prod_{j=-1}^{n-3}v_j\right)\ov{U}_{n,k}^{(1)}z_1\ov{U}_{n,k}^{(2)}
            z_2\ov{U}_{n,k}^{(3)}z_3\cdots,
\]
where $\bz := z_1z_2z_3\cdots$ is given by  $c_{\alpha_{n,1-k}}$ over the 
alphabet $\{(\ov{U}_{n,k-1})^{-1}, (U_{n-1,d_{n}-2})^{-1}\}$.
\end{lemma}
\begin{proof} Follows almost immediately from Proposition \ref{P:21.08.03(3)}. Indeed, 
$G^k(c_{\alpha_{n,1-k}}) = c_{\alpha_{n,1}}$ for $0 \leq k \leq d_{n+1}-1$, and hence  
\[
 c_\alpha = \left(\prod_{j=-1}^{n-3}v_j\right)c_{\alpha_{n,1}}\{v_{n-2}w_{n-1}, w_{n-2}v_{n-3}\}
                  =\left(\prod_{j=-1}^{n-3}v_j\right)c_{\alpha_{n,1-k}}\{v_{n-2}w_{n-1}, 
                   (v_{n-2}w_{n-1})^kw_{n-2}v_{n-3}\},
\]
where $v_{n-2}w_{n-1} = \ov{U}_{n,k}(\ov{U}_{n,k-1})^{-1}$ and 
\[
(v_{n-2}w_{n-1})^kw_{n-2}v_{n-3}  =   \ov{U}_{n,k}(v_{n-2})^{-1}w_{n-2}v_{n-3}  =   \ov{U}_{n,k}(U_{n-1,d_n-2})^{-1}.
\] 
\end{proof}

\section{Main result}

We are now equipped with the necessary tools to prove the main
result of this paper, which, in view of Corollary \ref{Cor:all_palindromes}, completely describes  occurrences of palindromes in $c_\alpha$. 

\begin{theorem} \label{T:21.08.03(1)}
Let $u$ be a palindromic factor of $c_\alpha$ with $|u|\geq2$. 
\begin{itemize}
\item[$(1)$] Suppose $u = vU_{n,k}\rev{v}$ for some $n\geq1$, where  
$v \suff v_{n-2}$, $v \ne v_{n-2}$, and $0 \leq k \leq d_{n+1}
-2$. Then 
\[
  c_\alpha = \left(\prod_{j=-1}^{n-2}v_j\right)v^{-1}u^{(1)}z_1u^{(2)}z_2
u^{(3)}z_3\cdots,
\]
where $\bz := z_1z_2z_3\cdots$ is given by 
$c_{\alpha_{n,-k}}$ over the alphabet 
$\{(vU_{n,{k-1}}\rev{v})^{-1}, \rev{v}^{-1}w_{n}{v}^{-1}\}$.
\item[$(2)$] Suppose $u = v\ov{U}_{n,k}\rev{v}$ for some $n\geq1$, where 
$v \suff w_{n-1}$, $v \ne w_{n-1}$, and $0 \leq k \leq d_{n+1}
-1$. Then   
\[
  c_\alpha = \left(\prod_{j=-1}^{n-3}v_j\right)v^{-1}u^{(1)}z_1u^{(2)}z_2
u^{(3)}z_3\cdots,
\]
where $\bz := z_1z_2z_3\cdots$ is given by 
$c_{\alpha_{n,1-k}}$ over the alphabet $\{(v\ov{U}_{n,k-1} \rev{v})^{-1}, (vU_{n-1,d_n-2}\rev{v})^{-1}\}$.
\end{itemize}
Moreover, if $u = a^k$ for some $k \in [2, d_1-1]$, then $c_\alpha = u^{(1)}z_1u^{(2)}z_2u^{(3)}z_3\cdots$, where $z_1z_2z_3\cdots$ is given by $c_{\alpha_{0,-k}}$ over the alphabet $\{(a^{k-1})^{-1},b\}$.
\end{theorem}
\begin{note} In regards to assertion (1), let us point out that $v$ is a suffix (and $\rev{v}$ is a prefix) of $w_{n}$ since $w_n = w_{n-2}v_{n-3}v_{n-2} = v_{n-2}v_{n-3}w_{n-2}$. Therefore, 
$\rev{v}^{-1}w_nv^{-1} \in \cAstar$ since $|v| < |v_{n-2}| = q_{n} -q_{n-1} \leq \frac{1}{2}q_n = \frac{1}{2}|w_n|$.  
\end{note}
\begin{proof}[Proof of Theorem $\ref{T:21.08.03(1)}$]
Assertions (1) and (2) are proved in a similar fashion, using Lemmas \ref{L:cool2} and \ref{L:cool3}  respectively, so we just give the proof of $(1)$. The last statement is trivial since $c_\alpha = G^k(c_{\alpha_{0,-k}}) = c_{\alpha_{0,-k}}\{a,a^kb\}$.

Suppose $u = vU_{n,k}\rev{v}$ for some $n\geq1$, where $v \suff v_{n-2}$, $v\ne v_{n-2}$ and $0 \leq k \leq d_{n+1}-2$. From Lemma \ref{L:cool2}, it follows that
\[
  c_\alpha = \left(\prod_{j=-1}^{n-2}v_j\right)v^{-1}(vU_{n,k}^{(1)}\rev{v})
\rev{v}^{-1}z_1v^{-1}(vU_{n,k}^{(2)}\rev{v})\rev{v}^{-1}z_2v^{-1}
(vU_{n,k}^{(3)}\rev{v})\rev{v}^{-1}z_3\cdots,
\]
where $\bz := z_1z_2z_3\cdots$ is given by $c_{\alpha_{n,-k}}$ over the 
alphabet $\{(U_{n,k-1})^{-1}, w_n\}$. Consequently,
since each occurrence of $u$ in $c_\alpha$ corresponds to an occurrence of $U_{n,k}$ in 
$c_\alpha$, we have 
\[
  c_\alpha = \left(\prod_{j=-1}^{n-2}v_j\right)v^{-1}u^{(1)}
\hat{z}_1u^{(2)}\hat{z}_2u^{(3)}\hat{z}_3\cdots,
\]
where $\hat{z}_{i} = \rev{v}^{-1}z_{i}v^{-1}, \mbox{for all}~ i\geq1$.
(Note that $\left(\prod_{j=-1}^{n-2}v_j\right)v^{-1} \in \cAstar$ since $v \suff v_{n-2}$.)
Thus, $\hat{\bz} := \hat{z}_1\hat{z}_2\hat{z}_3\cdots$ is given by $c_{\alpha_{n,-k}}$ over the alphabet 
$\{({v}U_{n,k-1}\rev{v})^{-1}, \rev{v}^{-1}w_nv^{-1}\}$. Indeed, $\hat{z}_i = \rev{v}^{-1}w_nv^{-1}$ if $z_i = w_n$, and  $\hat{z}_i = \rev{v}^{-1}(U_{n,k-1})^{-1}v^{-1} = ({v}U_{n,k-1}\rev{v})^{-1}$ if  
$z_i = (U_{n,k-1})^{-1}$.  
This completes the proof of (1). 
 
In part (2), note that $\left(\prod_{j=-1}^{n-3}v_j\right)v^{-1} \in \cAstar$, 
since $v$ is a proper suffix of $w_{n-1}$, and hence a suffix of $\prod_{j=-1}^{n-3}v_j = x^{-1}w_{n-1}$, 
where $x \in \cA$ (by Lemma \ref{L:2}$(5)$).

\end{proof}
\vspace{0.2cm}
\begin{example} \label{Ex:27.08.03}
Let us now demonstrate Theorem \ref{T:21.08.03(1)} for $c_\alpha$ with 
$\alpha = [0;\ov{2,1,3,1}] = (4\sqrt{5} - 5)/11$. In this case, we have 
\[
  c_\alpha = abaabaabaababaabaabaabaababaabaabaabaab\cdots.
\]
Also note that
\[
  w_{-1} = v_{-1} =  a, ~w_{0} = v_0 = b, ~w_{1} = aa, ~w_{2} = bab, 
  ~w_3 = aabaabaabaa, ~v_{1} = aabaabaa, ~v_{2} = bab.
\]
\begin{itemize}
\item[(i)] Consider the palindromic factor $u = baaw_2aab = baababaab$, where $v = baa \suff v_{1}$.  By Theorem \ref{T:21.08.03(1)}(1),
\begin{align*}
  c_\alpha &= v_{-1}v_0v_1(baa)^{-1}u^{(1)}z_1u^{(2)}z_2u^{(3)}z_3\cdots \\ 
           &= abaabaa(baababaab)z_1(baababaab)z_2(baababaab)z_3\cdots,
\end{align*}
where $\bz := z_1z_2z_3\cdots$ is given by $c_{\alpha_{3}}$ over the alphabet $\{aa,aabaa\}$. We have  $c_{\alpha_3} = c_{\alpha_{3,0}} = [0;\ov{1,2,1,3}] = \sqrt5 - 57/38$, and hence $c_{\alpha_3} = bbabbbabbba\cdots$. 
Thus, we can write
\[
  c_\alpha = abaabaa(baababaab)aabaa(baababaab)aabaa(baababaab)
           aa(baababaab)aabaa(baababaab)\cdots.
\]
\item[(ii)] Now consider the palindromic factor $u = U_{2,1} = (w_1v_{0})^1w_1 = aabaa$.   
By Theorem \ref{T:21.08.03(1)}(1), 
\[
  c_\alpha = v_{-1}v_0u^{(1)}z_1u^{(2)}z_2u^{(3)}z_3\cdots 
           = ab(aabaa)z_1(aabaa)z_2(aabaa)z_3\cdots,
\]
where $\bz := z_1z_2z_3\cdots$ is given by  
$c_{\alpha_{2,-1}}$ over the alphabet $\{(aa)^{-1},bab\}$. We have $\alpha_{2,-1} =
[0;2,\ov{1,2,1,3}] = (4\sqrt5 - 2)/{19}$, and therefore
\[
  c_{\alpha_{2,-1}} = abaabaababaabaabaababaabaabaa\cdots.
\]
Hence, we can write 
\begin{align*}
c_\alpha =
&~ab(aab\underline{aa)baa}bab(aab\underline{aa)b(aa}baa)bab(aa
               b\underline{aa)b(aa}baa)bab(aab\underline{aa)baa} \\
            &~bab(aab\underline{aa)b(aa}baa)bab(aab\underline{aa)baa}\cdots.
\end{align*}
Notice that $u$ has a unique overlap factor $aa$.
\item[(iii)] Let us now consider the palindromic factor $u = aU_{2,2}a = a(v_{0}w_1)^2v_0a = abaabaaba$, where $a \suff w_{1} =aa$. Observe that 
\[
 (vU_{2,1}\rev{v})^{-1} = (a(baa)^1ba)^{-1} = (abaaba)^{-1} 
 \quad \mbox{and} \quad (vU_{1,d_2 -2}\rev{v})^{-1} = (a(ba)^{-1}ba) = a^{-1}.
\]
Thus, by Theorem \ref{T:21.08.03(1)}(2), we have
\[
  c_\alpha = v_{-1}a^{-1}u^{(1)}z_1u^{(2)}z_2u^{(3)}z_3\cdots
           = (abaabaaba)z_1(abaabaaba)z_2(abaabaaba)z_3\cdots,
\]
where $\bz := z_1z_2z_3\cdots$ is given by 
$c_{\alpha_{2,-1}}$ over the alphabet $\{(abaaba)^{-1},a^{-1}\}$; whence 
\begin{align*}
c_\alpha =
&~(aba\underline{abaaba)ab(a}ba\underline{aba[aba)aba}ab
           \underline{a]ba(aba[aba}aba)ab\underline{a]ba(abaaba}a \\
            &~b\underline{a)ba(aba[aba}aba)ab\underline{a]baabaaba}\cdots.
\end{align*}
In this case, $u$ has two overlap factors: $abaaba$ and $a$.
\end{itemize} \qed
\end{example}

Let us now denote by occ$_i(u)$ the position of the $i$-th
occurrence of $u$ in $c_\alpha$, i.e., if $c_\alpha = zu\bx$ for some 
$z \in \cAstar$, $\bx \in \cAw$ such that $|zu|_u = i$, then
occ$_i(u) = |z|$. With this notation, a given factor $u$ of $c_\alpha$ occurs at 
precisely the positions (occ$_i(u))_{i \geq1}$ in $c_\alpha$. 

The following corollary of Theorem \ref{T:21.08.03(1)} gives the 
exact positions at which palindromes occur in $c_\alpha$.

\begin{corollary} \label{Cor:25.07.03(2)}
Let $u$ be a palindromic factor of $c_\alpha$ with $|u|\geq2$.
\begin{itemize}
\item[$(1)$] Suppose $u = vU_{n,k}\rev{v}$ for some $n\geq1$, where 
$v \suff v_{n-2}$, $v \ne v_{n-2}$, and $0 \leq k \leq d_{n+1} - 
2$. Then \emph{occ}$_1(u) = \frac{1}{2}((k+2)q_{n} + q_{n-1}-|u|-2)$  and, for  all $i\geq1$, 
\[
 \mbox{\emph{occ}$_{i+1}(u) =$ \emph{occ}$_{i}(u) + P_i$},
\] 
where $(P_i)_{i\geq1}$  is given by   
$c_{\alpha_{n,-k}}$ over the alphabet $\{q_{n},(k+1)q_{n} + q_{n-1}\}$.
\item[$(2)$] Suppose $u = vU_{n,k}\rev{v}$ for some $n\geq1$, where 
$v \suff w_{n-1}$, $v \ne w_{n-1}$, and $0 \leq k \leq d_{n+1} -1$. Then \emph{occ}$_1(u) = \frac{1}{2}((k+1)q_{n} + q_{n-1}-|u|-2)$ and, for all $i\geq1$,  
\[
 \mbox{\emph{occ}$_{i+1}(u) =$ \emph{occ}$_{i}(u) + P_i$}, 
\] 
where $(P_i)_{i\geq1}$  is given by $c_{\alpha_{n,1-k}}$ over 
the alphabet $\{q_{n}, kq_{n} + q_{n-1}\}$.
\end{itemize}
Moreover, if $u = a^k$ for some $k \in [2,d_1-1]$, then \emph{occ}$_1(u) = 0$ and \emph{occ}$_{i+1}(u) =$ \emph{occ}$_{i}(u) + P_i$ for all $i\geq1$, where $(P_i)_{i\geq1}$ 
is given by $c_{\alpha_{0,-k}}$ over the alphabet $\{1,k+1\}$.
\end{corollary}
\begin{proof}
As with Theorem \ref{T:21.08.03(1)}, the proofs of (1) and (2) are much the same, so we just give the proof of $(1)$. The proof of the last statement is trivial. 

Suppose $u = vU_{n,k}\rev{v}$ for some $n\geq1$, where $v \suff v_{n-2}$, $v\ne v_{n-2}$ and $0\leq k\leq d_{n+1}-2$.  Theorem \ref{T:21.08.03(1)}$(1)$ shows that
\begin{equation} \label{eq:dododo}
  c_\alpha = \left(\prod_{j=-1}^{n-2}v_j\right)v^{-1}u^{(1)}z_1u^{(2)}z_2
u^{(3)}z_3\cdots,
\end{equation}
where $\bz := z_1z_2z_3\cdots$ is given by $c_{\alpha_{n,-k}}$ over the alphabet 
$\{(vU_{n,k-1}\rev{v})^{-1}, \rev{v}^{-1}w_nv^{-1}\}$. 
Observe that
\[
|U_{n,k}| = (k+1)|w_{n-1}| + k|v_{n-2}| = (k+1)q_{n-1} + k(q_n - q_{n-1}) = k q_n + q_{n-1}, 
\]
and hence
\[
 |v| = \frac{1}{2}(|u| - |U_{n,k}|) = \frac{1}{2}(|u| - kq_{n} - q_{n-1}).
\] 
Also recall that if $x$ is the first letter of $w_{n}$, then $x^{-1}w_{n} = \prod_{j=-1}^{n-2}v_j$.
Therefore, since $v$ is a proper suffix of $v_{n-2}$, we have
\[
  \left|\left(\prod_{j=-1}^{n-2}v_j\right)v^{-1}\right| = |x^{-1}w_{n}| - |v| = q_{n} - 1- 
\frac{1}{2}(|u| - kq_{n} - q_{n-1}).
\]
Hence, the first occurrence of $u$ in $c_\alpha$ is at position 
\[
  \textrm{occ}_1(u) = \frac{1}{2}((k+2)q_{n} + q_{n-1} - |u| - 2).
\]
Furthermore, 
\[
|{v}U_{n,k-1}\rev{v}| = (k-1)q_{n} + q_{n-1} + 2|v|
                          = (k-1)q_{n} + q_{n-1} + (|u| - kq_{n} - q_{n-1})
                          = |u| - q_{n} 
\] 
and 
\[
|\rev{v}^{-1}w_nv^{-1}| = q_n - 2|v| = q_n - (|u| - kq_{n} -q_{n-1}) = (k+1)q_{n} + q_{n-1}-|u|.
\]
Thus, it follows from \eqref{eq:dododo} that  occ$_{i+1}(u) =$ occ$_{i}(u) + P_i$ for all $i\geq1$,  where 
$(P_i)_{i\geq1}$ is the characteristic Sturmian word of slope
$\alpha_{n,-k}$ over the alphabet $\{q_{n}, (k+1)q_n + q_{n-1}\}$. 
\end{proof}

\begin{example} Let $\alpha = [0;\ov{2,1,3,1}] = (4\sqrt5 - 5)/11$ and consider the 
palindromic factor $u$ of $c_\alpha$ given by $u = U_{2,1} = aabaa$. According to Corollary 
\ref{Cor:25.07.03(2)}, one should find that $u$ first occurs at position
\[
  \textrm{occ}_1(u) = \frac{1}{2}(3q_2 - q_1 - 5 - 2) = \frac{1}{2}(9+2-7) = 2,
\]
followed by the positions occ$_{i+1}(u) =$ occ$_i(u) + P_i$ for each $i\geq1$, 
where $(P_i)_{i\geq1}$ is the characteristic Sturmian word of slope $\alpha_{2,-1} =
[0;2,\ov{1,2,1,3}]$ over the alphabet
$\{q_2,2q_2+q_1\} = \{3,8\}$; that is, $(P_i)_{i\geq1} =
(3,8,3,3,8,3,3,8,3,8,3,3,8,3,3,8,\ldots)$. Indeed, from Example
\ref{Ex:27.08.03}$(2)$, we have
\begin{align*}
c_\alpha =
&~ab(aab\underline{aa)baa}bab(aab\underline{aa)b(aa}baa)bab(aa
               b\underline{aa)b(aa}baa)bab(aab\underline{aa)baa} \\
            &~bab(aab\underline{aa)b(aa}baa)bab(aab\underline{aa)baa}\cdots, 
\end{align*}
from which it is evident that $u = aabaa$ occurs at positions
  2, 5, 13, 16, 19, 27, 30, 33, 41, 44, 52, 55, 58, 66, $\ldots$ .
\end{example}

\begin{remark}
In general, if $u_1$ and $u_2$ are the two return words of a factor $u$ of $c_\alpha$, it is clear that 
\[
  \mbox{occ$_{i+1}(u) =$ {occ}$_i(u) + |u_{j_i}|$ where $j_i = 1$ or $2$}.
\] 
In particular, the sequence $(j_i)_{i\geq1}$ is a Sturmian word over 
the alphabet $\{1,2\}$ (see \cite{iFlV02gene} or Section \ref{S:return}).  In the case when 
$u$ is a palindromic factor of $c_\alpha$, Corollary 
\ref{Cor:25.07.03(2)} shows that the sequence $(j_{i})_{i\geq1}$ is given by 
$c_{\alpha_{n,-k}}$ over the alphabet $\{1,2\}$, for some $n\in\NN$ and $0 \leq k \leq d_{n+1}-1$.  For example, the two return words of $w_n$ ($= U_{n+1,0}$) are 
$u_1 = w_nv_{n-1}$ and $u_2 = w_{n}w_{n+1}$, where 
\[
  |u_1| = q_{n+1} \quad \mbox{and} \quad |u_2| = q_{n+1} + q_n. 
\]  
From Corollary \ref{Cor:25.07.03(2)}, occ$_1(w_n) = q_{n+1}-1$ and $\mbox{occ}_{i+1}(w_n) = 
\mbox{occ}_i(w_n) + |u_{j_i}|$ for each $i\geq1$, where $(j_i)_{i\geq1}$ is given by $c_{\alpha_{n+1}}$ over the alphabet $\{1,2\}$.
\end{remark}

\section{Occurrences of factors of length $q_n$ in $c_\alpha$}

In this last section, we determine the structure of all factors of length
$q_n$ of $c_\alpha$ with respect to the singular words $w_n$, $w_{n-1}$, and
$v_{n-2}$. Subsequently, using some results from Section 4, we 
completely describe where factors of length $q_n$ occur in
$c_\alpha$.

Let $w = x_1x_2\cdots x_m \in \cAstar$, each $x_i \in \cA$, and
let $k \in \NN$ with $0 \leq k \leq m - 1$. The \emph{$k$-{th}
conjugate} of $w$ is the word  
$C_k(w) := x_{k+1}x_{k+2}\cdots x_{m}x_1x_2\cdots x_k$.  
Further, we conventionally set $C_{-k}(w) = C_{|w| - k}(w)$ and define $C(w) := \{C_k(w) : ~0 \leq k \leq |w| - 1 \}$. 

One can easily prove that any conjugate of $s_n$ is a factor of
$c_\alpha$. Certainly, $C(s_{-1}) = \{b\}$ and, for $n \geq 0$, 
\[
  s_{n+3} = s_{n+2}^{d_{n+3}}s_{n+1}
          = (s_{n+1}^{d_{n+2}}s_n)^{d_{n+3}}s_n^{d_{n+1}}s_{n-1}
          = (s_{n+1}^{d_{n+2}}s_n)^{d_{n+3}-1}s_{n+1}^{d_{n+2}}
              s_{n}^{d_{n+1}+1}s_{n-1},
\]
where $d_{n+1} + 1 \geq 2$. Thus, $s_n^2$ is a factor of $c_\alpha$, and hence the 
claim is proved since any conjugate of $s_n$ is a factor of $s_n^2$. 

Now, each $s_n$ is a primitive word \cite{aDfM94some}, i.e., $s_n$
cannot be written as a non-trivial integer power of a shorter word.
Consequently, $s_n$ has $q_n$ distinct conjugates, i.e., $|C(s_n)|
= q_n$. Furthermore, from the above observation, $C(s_n)$ is a set
of factors of $c_\alpha$. It is therefore deduced that the set of all
factors of length $q_n$ of $c_\alpha$ consists of $C(s_n)$ and
$w_n$. That is,
\[
  \Omega_{q_n}(c_\alpha) = C(s_n) \cup \{w_n\}.
\]
Indeed, since $c_\alpha$ is a Sturmian word, it must have
exactly $q_n + 1$ distinct factors of length $q_n$. 

\begin{lemma} \label{L:5} For any $n \geq 1$, 
$C_{q_n-1}(s_n) = w_{n-1}v_{n-2}$ and $C_{q_{n-1}-1}(s_n) = v_{n-2}w_{n-1}$. 
Moreover, 
\begin{itemize}
\item[$(1)$] for $0 \leq k \leq q_{n-1} - 2$, $C_{k}(s_n) = uv_{n-2}v$, where 
$vu = w_{n-1}$ and $|v| = k+1;$
\item[$(2)$] for $q_{n-1} - 1 \leq k \leq
q_{n} - 1$, $C_{k}(s_n) = uw_{n-1}v$, where $vu = v_{n-2}$ and $|v|  = k + 1 - q_{n-1}$.
\end{itemize}
\end{lemma}
\begin{proof}
By Lemma \ref{L:2}$(1)$, $C_{-1}(s_n) = C_{q_n - 1}(s_n) =
w_{n-1}v_{n-2}$. Therefore, 
$C_{q_{n-1}-1}(s_n) = v_{n-2}w_{n-1}$  since $|w_{n-1}| = q_{n-1}$. Assertions (1) and (2)  
follow immediately.
\end{proof}
Accordingly, a factor of length $q_n$ of $c_\alpha$ is either
$w_n$, or has at least one of the words $v_{n-2}$ and $w_{n-1}$ as a 
factor. We shall now establish two different decompositions of $c_\alpha$,
which show exactly where conjugates of $s_n$ occur in $c_\alpha$.

\begin{theorem} \label{T:09.09.03} Let $n\geq 1$.
\begin{itemize}
\item[$(1)$] Suppose $w = C_k(s_n)$ for some $k \in [0, q_{n-1}-2]$, so that $w = uv_{n-2}v$,  where $vu = w_{n-1}$ and $|v| = k +1$. Then 
\[
  c_\alpha = \left(\prod_{j = -1}^{n-3}v_j\right)u^{-1}w^{(1)}z_1w^{(2)}z_2
              w^{(3)}z_3\cdots,
\]
where $\bz := z_1z_2z_3\cdots$ is given by $c_{\alpha_{n,1}}$ over the alphabet $\{\empt,
(uU_{n-1,d_n-2}v)^{-1}\}$.
\item[$(2)$] Suppose $w = C_k(s_n)$ for some $k \in [q_{n-1}-1, q_{n}-1]$, so that $w = uw_{n-1}v$, where $vu = v_{n-2}$ and $|v| = k+1- q_{n-1} $. Then 
\[
  c_\alpha = \left(\prod_{j = -1}^{n-2}v_j\right)u^{-1}w^{(1)}z_1w^{(2)}z_2
              w^{(3)}z_3\cdots,
\]
where $\bz := z_1z_2z_3\cdots$ is given by $c_{\alpha_{n}}$ over the alphabet 
$\{\empt, v^{-1}w_{n}u^{-1}\}$.
\end{itemize}
\end{theorem}
\begin{proof}
Using decompositions \eqref{eq:v_{n-1}} and \eqref{eq:w_n} (consequences of 
Lemma \ref{L:cool2}), the proof follows along exactly the same lines as the proof of 
Theorem \ref{T:21.08.03(1)}.
\end{proof}

In light of Theorem \ref{T:09.09.03} and the $w_n$-decomposition of 
$c_\alpha$ given by \eqref{eq:w_n} (together with the fact that $\Omega_{q_n}(c_\alpha) = C(s_n) \cup \{w_n\}$), we have now shown precisely where each factor of length $q_n$ occurs in $c_\alpha$. It is important to note that it follows from
Proposition \ref{P:bCzW03-10} that a factor $w$ of length $q_n$ does 
not have overlap in $c_\alpha$ if and only if $w = w_{n}$, or $w = C_k(s_n)$ for some $k \in [q_{n-1}-1, q_n -1]$. Certainly, if $w$ takes the latter form, then $w = uw_{n-1}v$ with $vu=v_{n-2}$ and 
$|v| = k+1- q_{n-1}$. In this case, Theorem 
\ref{T:09.09.03}(2) shows that $w$ does not have overlapping occurrences since $w_n = (vu)v_{n-3}w_{n-2} = w_{n-2}v_{n-3}(vu)$, where $vu = v_{n-2}$, and hence $v^{-1}w_nu^{-1} \in \cAstar$.

\begin{example} \label{Ex:09.09.03} 
Suppose $\alpha = [0;\ov{2,1}] = (\sqrt{3} - 1)/{2}$. Then
\[
  c_\alpha = abaabaababaabaabaababaabaabaababaabaababaabaabaababaabaabaababa
             \cdots.
\]
Let us demonstrate the above theorem by considering the first two conjugates of $s_3 = abaabaab$; 
namely, $C_1(s_3)$ ($=C_{q_2-2}(s_3)$) and $C_2(s_3)$ ($=C_{q_2-1}(s_3)$).  
First observe that 
\[ 
w_{-1} = v_{-1} = a, ~w_{0} = v_0 =  b, ~w_{1} = aa, ~w_{2} = bab, 
 ~w_3 = aabaabaa, ~v_{1} = aabaa, ~v_{2} = bab.
\]
\begin{itemize}
\item[$(1)$] Let $w = C_1(s_3)$; the first conjugate of $s_3$. We have 
$w = baabaaba = uv_1v$, where $u = b$, $v=ba$ and $vu = bab = w_2$. 
Hence, by Theorem \ref{T:09.09.03}$(1)$,
\[
  c_\alpha = v_{-1}v_0b^{-1}w^{(1)}z_1w^{(2)}z_2w^{(3)}z_3\cdots 
           = a(baabaaba)z_1(baabaaba)z_2(baabaaba)z_3\cdots,
\]
where $\bz := z_1z_2z_3\cdots$ is given by $c_{\alpha_{3,1}}$ over the alphabet $\{\empt,(baaba)^{-1}\}$. (Note that $(uU_{2,d_3-2}v)^{-1} = (b(aab)^0aaba)^{-1} = (baaba)^{-1}$.)  Since $\alpha_{3,1} =
[0;2,\ov{2,1}] = 1 - \sqrt{3}/3$, we have $c_{\alpha_{3,1}} = ababaabababaab\cdots$, and thus 
we can write 
\begin{align*}
  c_\alpha = &~a(baabaaba)(baa\underline{baaba)aba}(baa\underline{baaba)aba}
           (baabaaba)(baa\underline{baaba)aba} \\
             &~(baa\underline{baaba)aba}(baa\underline{baaba)aba}(baabaaba)
           (baa\underline{baaba)aba}\cdots.
\end{align*}
\item[$(2)$] Now let $w = C_2(s_3)$; the second conjugate of $s_3$. Then $w = aabaabab = v_1w_2$, and it follows from Theorem \ref{T:09.09.03}(2) that
\[
  c_\alpha = v_{-1}v_0v_1{v_1}^{-1}w^{(1)}z_1w^{(2)}z_2w^{(3)}z_3\cdots 
           = ab(aabaabab)z_1(aabaabab)z_2(aabaabab)z_3\cdots,
\]
where $\bz := z_1z_2z_3\cdots$ is given by $c_{\alpha_{3}}$ over the alphabet 
$\{\empt, w_3(aabaa)^{-1}\} = \{\empt, aab\}$. Note that  
$\alpha_{3} = [0;\ov{1,2}] = \sqrt{3} - 1$, and hence 
$c_{\alpha_{3}} = bbabbbabbba\cdots$.  Therefore,  
\begin{align*}
  c_\alpha = &~ab(aabaabab)aab(aabaabab)aab(aabaabab)(aabaabab)aab(aa
             baabab) \\
             &~aab(aabaabab)aab(aabaabab)(aabaabab)aab(aabaabab)aab\cdots.
\end{align*}
\end{itemize} 
\end{example}

\begin{remark}
From Theorem \ref{T:09.09.03}, one can easily deduce Lemma \ref{L:21.08.03(2)} (i.e., 
$c_\alpha = c_{\alpha_{n,1}}\{s_n,s_{n-1}\}$ for all $n\geq 1$), as follows. Observe that $s_n = C_0(s_n) = uv_{n-2}v$, where $u = y^{-1}w_{n-1}$ 
and $v = y$. Similarly, $s_{n-1} = x^{-1}w_{n-2}v_{n-3}x$, where $x \in \cA$
($x\ne y$). Hence, 
\[
  c_\alpha = \left(\prod_{j=-1}^{n-3}v_j\right)(y^{-1}w_{n-1})^{-1}
             s_{n}^{(1)}z_1s_{n}^{(2)}z_2s_{n}^{(3)}z_3\cdots
           = s_{n}^{(1)}z_1s_{n}^{(2)}z_2s_{n}^{(3)}z_3\cdots,
\]
where $\bz := z_1z_2z_3\cdots$ is given by $c_{\alpha_{n,1}}$ over the alphabet 
$\{\empt,
(y^{-1}w_{n-1}U_{n-1,d_n-2}y)^{-1}\}$. That is,
\[
  c_\alpha = c_{\alpha_{n,1}}\{s_n, s_n(y^{-1}w_{n-1}(w_{n-2}v_{n-3})^{-1}
  v_{n-2}y)^{-1}\}.
\]
Using Lemma \ref{L:2}, we have 
\begin{align*}
  s_n(y^{-1}w_{n-1}(w_{n-2}v_{n-3})^{-1}v_{n-2}y)^{-1} &= s_n(s_{n-1}x^{-1}(w_{n-2}v_{n-3})^{-1}v_{n-2}y)^{-1} \\
  &= s_ny^{-1}(v_{n-2})^{-1}(w_{n-2}v_{n-3})x(s_{n-1})^{-1} \\
  &= x^{-1}w_n(v_{n-2})^{-1}xs_{n-1}(s_{n-1})^{-1} \\ 
  &=x^{-1}w_{n-2}v_{n-3}x \\
  &=s_{n-1}, 
\end{align*}
and therefore $c_\alpha = c_{\alpha_{n,1}}\{s_{n}, s_{n-1}\}$, as required.
\end{remark}

We finish with a corollary of Theorem \ref{T:09.09.03} (\emph{cf.} Corollary
\ref{Cor:25.07.03(2)}).

\begin{corollary} \label{Cor:20.03.04}
Let $n \geq1$ and suppose $w = C_k(s_n)$ for some $k \in [0, q_{n}-1]$. Then \emph{occ}$_1(w) = k$ and, for all $i\geq1$, \emph{occ}$_{i+1}(w) =$ \emph{occ}$_{i}(w) + P_i$,  
where $(P_i)_{i\geq1}$  is given by: 
\begin{itemize}
\item $c_{\alpha_{n,1}}$ over the alphabet $\{q_n, q_{n-1}\}$ if $k \in[0, q_{n-1}-2],$ or 
\item  $c_{\alpha_{n}}$ over the alphabet $\{q_n, q_n + q_{n-1}\}$ if $k \in [q_{n-1}-1, q_n-1]$.
\end{itemize} \qed
\end{corollary}

\section{Acknowledgements} Special thanks to Bob Clarke and Alison Wolff for
their support and encouragement. Thanks also to the two anonymous referees for their helpful suggestions and comments. This research was supported by the George Fraser Scholarship of The University of Adelaide.

%------------------References----------------------------

\end{document}